\documentclass[amssymb,amstex,amsmath,10pt,a4paper]{amsart}
\usepackage{amsmath,amssymb,latexsym,amsfonts,url}
\usepackage{graphicx}

\newtheorem{thm}{Theorem}[section]

\newtheorem{lem}[thm]{Lemma}

\newcommand{\Div}{{\rm div}}

\newenvironment{Rmk}{\noindent {\it Remark.}}{\\}
\numberwithin{equation}{section}

\begin{document}
\title[Sphere-foliated hypersurfaces in product spaces]{Sphere-foliated minimal and constant mean curvature hypersurfaces in product spaces}
\author{Keomkyo Seo}

\begin{abstract}
In this paper, we prove that minimal hypersurfaces when $n\geq 3$ and nonzero constant mean curvature hypersurfaces when $n\geq2$ foliated by spheres in parallel horizontal hyperplanes in ${\mathbb{H}}^n \times \mathbb{R}$ must be rotationally symmetric. \\

\end{abstract}
\maketitle
\section{Introduction}
In $\mathbb{R}^3$, catenoids and Riemann's examples are the only
minimal surfaces foliated by circles. However, in higher-dimensional
Euclidean space, there are no examples of non-rotationally symmetric
minimal hypersurfaces such as Riemann's examples in $\mathbb{R}^3$.
In 1991, Jagy \cite{Jagy91} proved that
if $M^n$ is a minimal hypersurface in $\mathbb{R}^{n+1}$ $(n \geq
3 )$ and foliated by $(n-1)$-dimensional spheres in parallel
hyperplanes, then $M^n$ is rotationally symmetric about the axis
containing the centers of all the spheres. This result has been
generalized to other spaces forms: sphere, the hyperbolic and
Lorentz-Minkowski space.(See \cite{Jagy98}, \cite{Lopez99}, and
\cite{Park02}.)

In $\mathbb{H}^2 \times \mathbb{R}$, Nelli and
Rosenberg \cite{NR02} found a rotationally symmetric minimal
surface which is called a $catenoid$ in $\mathbb{H}^2 \times
\mathbb{R}$. In \cite{Hauswirth06}, Hauswirth provided several
examples of minimal surfaces foliated by horizontal curves of
constant curvature in $\mathbb{H}^2 \times \mathbb{R}$. In
particular, he constructed a two-parameter family of Riemann type
surfaces. Recently, B\'{e}rard and Sa Earp \cite{BSE08} obtained
some results on total curvature and index of higher-dimensional catenoids in
${\mathbb{H}}^n \times \mathbb{R}$. On the other hand, Nelli et
al. \cite{Nelli et al08} described the geometric behavior of
rotationally symmetric constant mean curvature surfaces in
${\mathbb{H}}^2 \times \mathbb{R}$. They showed that for
$|H|>1/2$, the properties of rotationally symmetric constant mean
curvature surfaces in ${\mathbb{H}}^2 \times \mathbb{R}$ are
analogous to those of the Delaunay surfaces in $\mathbb{R}^3$.
Rotationally symmetric constant mean curvature surfaces in
${\mathbb{H}}^2 \times \mathbb{R}$ have been studied in
\cite{AR04, CdL07, Hsiang89, SET05}.

Throughout this paper, we consider the upper half-space model of hyperbolic space
\begin{eqnarray*}
\mathbb{H}^{n} = \{(x_1, \cdots, x_{n})\in \mathbb{R}^{n} : x_{n} > 0 \}
\end{eqnarray*}
equipped with the metric $\displaystyle{ds^2 = \frac{(dx_1)^2 + \cdots + (dx_{n})^2}{{x_n}^2}}$.
For a product space $\mathbb{H}^n \times \mathbb{R}$, we fix the metric
$ds^2 +\varepsilon dt^2 \, (\varepsilon=\pm1)$. This metric is called $Riemannian$ if $\varepsilon=1$ and
$Lorentzian$ if $\varepsilon=-1$.

In this paper, we study hypersurfaces foliated by
$(n-1)$-dimensional spheres lying in parallel hyperplanes in some
Riemannian and Lorentzian product spaces. In Section 2, we shall
prove that minimal hypersurfaces $(n \geq 3)$ and non-zero
constant mean curvature hypersurfaces $(n \geq 2)$ foliated by
$(n-1)$-dimensional spheres in parallel horizontal hyperplanes in the
Riemannian product ${\mathbb{H}}^n \times \mathbb{R}$ should be rotationally
symmetric.(Theorem \ref{thm:foliation}) As a consquence, one can see that there
is no Riemann type minimal hypersurfaces foliated by
$(n-1)$-dimensional spheres in ${\mathbb{H}}^n \times \mathbb{R}$
for $n\geq 3$. We shall use Jagy's idea \cite{Jagy91} to prove this
result. (See also \cite{Lopez99}.) A key ingredient of the proof
is the following. We describe a hypersurface $M$ in
${\mathbb{H}}^n \times \mathbb{R}$ locally as the level set for a smooth
function $f$. If we orient $M$ by the unit normal vector field
$\displaystyle{N=-\frac{\nabla f}{|\nabla f|}}$, then the mean curvature $H$ is
given by
\begin{eqnarray}\label{eqn:H}
nH = -\Div \frac{\nabla f}{|\nabla f|},
\end{eqnarray}
where $\nabla$ and $\Div$ denote the gradient and divergence in ${\mathbb{H}}^n \times
\mathbb{R}$, respectively. A straightforward computation using the fact that $M$ is
foliated by spheres in parallel horizontal hyperplanes gives us the
conclusion. In Section 3, applying the similar arguments as in Section 2,
we prove an analogue in the Lorentzian product ${\mathbb{H}}^n \times \mathbb{R}$.

The author would like to thank the referee for useful suggestions on improving the presentation of this paper.

\section{Sphere-foliated hypersurfaces in the Riemannian product $\mathbb{H}^n \times
\mathbb{R}$}

In $\mathbb{H}^n \times \mathbb{R}$, a one-parameter family of hyperplanes $\mathbb{H}^n \times
\{t\}$ for $t\in \mathbb{R}$ are called $parallel$ $horizontal$ hyperplanes.
We will deal with hypersurfaces foliated by spheres in parallel horizontal hyperplanes in the
Riemannian product space $\mathbb{H}^n \times \mathbb{R}$.

\begin{thm}\label{thm:foliation}
Let $M$ be an $n$-dimensional hypersurface with constant mean
curvature $H$ in the Riemannian product $\mathbb{H}^n \times
\mathbb{R}$ and foliated by spheres in parallel horizontal
hyperplanes. If $H \neq 0$ or $H=0$ and $n \geq 3$, then $M$ is a
rotationally symmetric hypersurface.
\end{thm}

Before proving the above theorem, we need the following well-known fact.

\begin{lem}[\cite{Stahl}, p.81-82] \label{lem:center}
If an $(n-1)$-dimensional sphere has Euclidean center $(0,\ldots, 0, k) \in \mathbb{R}^{n}_+ :=\{(x_1, \cdots, x_{n}):x_n >0 \}$ and a Euclidean radius $r$,
then it has the hyperbolic center $(0,\ldots, 0, K) \in \mathbb{R}^n_+$ and the hyperbolic radius $R$, where
\begin{eqnarray*}
 K=\sqrt{k^2-r^2} \, \, \mbox{and } \,\, R=\frac{1}{2}\ln \frac{k+r}{k-r}.
\end{eqnarray*}
\end{lem}
{\bf Proof of Theorem \ref{thm:foliation}.}
Let $P_{t_1} = \mathbb{H}^n \times \{t_1\}$ and $P_{t_2} =
\mathbb{H}^n \times \{t_2\}$ be two horizontal hyperplanes of the
foliation for $t_1 < t_2$. Let $M^*$ be the piece of $M$ between $P_{t_1}$ and $P_{t_2}$. The boundary $\partial M^*$ of $M^*$ consists of two $(n-1)$-dimensional spheres $(M^*\cap P_{t_1}) \cup (M^*\cap P_{t_2})$. After an isometric transformation in $\mathbb{H}^n \times \mathbb{R}$, we may assume that the hyperbolic centers of the two boundary spheres are given by $(0, \cdots, 0, k_1, t_1)$ and $(0, \cdots, 0, k_2, t_2)$ in $\mathbb{R}^{n}_+ \times \mathbb{R}$ for some constants $k_1, k_2 >0$, respectively. Note that these two boundary spheres are symmetric to the hyperplanes $\{x_1=0\}, \cdots, \{x_{n-1}=0\}$. The well-known Aleksandrov reflection principle shows that $M^*$ inherits the symmetries of its boundary $\partial M^*$. Therefore, for each $t_1 \leq t \leq t_2$, the hyperbolic center of each level $M \cap \{x_{n+1}=t\}$ is symmetric to the hyperplanes $\{x_1=0\}, \cdots, \{x_{n-1}=0\}$. Hence it follows that the hyperbolic center of each level lies in the $2$-plane $\{x_1= \cdots=x_{n-1}=0\}$.


Using Lemma \ref{lem:center}, we parametrize the hyperbolic centers of the
$(n-1)$-dimensional spheres by $t \longmapsto (0, \ldots, 0, K(t),
t) \in \mathbb{H}^n \times
\mathbb{R}$ for $t\in [t_1, t_2]$, and hence the Euclidean centers of the
spheres by $t \longmapsto (0, \ldots, 0, k(t), t)\in \mathbb{H}^n \times
\mathbb{R}$. Then it
follows from Lemma \ref{lem:center} that
\begin{eqnarray*}
K(t) = \sqrt{k(t)^2 - r(t)^2},
\end{eqnarray*}
where $r(t)$ is the Euclidean radius of $M\cap \mathbb{H}^n \times
\{t\}$. Note that $K(t)>0$ for $t\in [t_1, t_2]$. Moreover, we see
 that $M^*$ is the level set $\{f=0\}$ of a function $f$ given by
\begin{eqnarray} \label{eqn:function}
f(x_1, \ldots, x_n, t) = \sum_{i=1}^{n-1} x_i^2 + (x_n -k(t))^2
-r(t)^2.
\end{eqnarray}
To prove the theorem, it is sufficient to show that
\begin{eqnarray*}
\frac{dK(t)}{dt}=0,
\end{eqnarray*}
 which means that $M$ is a rotationally
symmetric hypersurface whose rotation axis is the
geodesic $\gamma (t) = \{(0, \ldots, 0, K, t)\}$ for some constant
$K$. Note that the metric on $\mathbb{H}^n  \times \mathbb{R}$ is given by
\begin{eqnarray*}
\sum_{i,j} g_{ij} dx_i \otimes dx_j = \frac{1}{x_n^2}dx_1^2 + \cdots + \frac{1}{x_n^2}dx_n^2+dt^2 .
\end{eqnarray*}
Since
\begin{eqnarray*}
\nabla f &=& \sum_{i,j} g^{ij}\frac{\partial f}{\partial x_i}\frac{\partial}{\partial x_j} \\
 &=& (2x_n^2x_1, \cdots, 2x_n^2 x_{n-1}, 2x_n^2 (x_n - k), -2(x_n -k)k' - 2r r'),
\end{eqnarray*}
we have
\begin{eqnarray*}
|\nabla f|^2 = 4\Big(x_n^2x_1^2+ \cdots+ x_n^2 x_{n-1}^2+ x_n^2 (x_n - k)^2 + ((x_n -k)k' + r r')^2\Big) .
\end{eqnarray*}
Now we compute the mean curvature of $M^*$ using the equation (\ref{eqn:H}).
\begin{eqnarray*}
-nH = \Div \frac{\nabla f}{|\nabla f|} &=& \sum_{i,j} \frac{1}{\sqrt{g}} \frac{\partial}{\partial x_j}\Big(\sqrt{g}\frac{g^{ij}}{|\nabla f|}\frac{\partial f}{\partial x_i}\Big)\\
&=& \sum_{j} \frac{\partial Z^j}{\partial x_j} + \sum_{j} \frac{1}{\sqrt{g}}\Big(\frac{\partial \sqrt{g}}{\partial x_j}\Big)Z^j,
\end{eqnarray*}
where $\displaystyle{Z^j = \sum_{i} \frac{g^{ij}}{|\nabla f|}\frac{\partial f}{\partial x_i}}$ and $g=\det(g_{ij})=x_n^{-2n}$. Then we have
\begin{align*}
-nH &= \frac{x_n^2}{S} - \frac{(x_n^2 x_1)^2}{S^3}+ \cdots + \frac{x_n^2}{S} - \frac{(x_n^2 x_{n-1})^2}{S^3} \\
& \quad + \frac{x_n^2+2x_n (x_n - k)}{S} -
\frac{x_n^2 (x_n - k)\{x_n r^2 + x_n^2 (x_n - k)+Ak'\}}{S^3}\\
& \quad + \frac{B}{S} - \frac{A\{x_n^2(x_n - k)k'+ AB\}}{S^3}\\
& \quad + x_n^n (-n x_n^{-n-1}) \frac{x_n^2 (x_n - k)}{S},
\end{align*}
where $A=(x_n -k)k' + r r',$ $B=k'^2-(x_n - k)k'' - r'^2 - rr'',$ and\\
\begin{eqnarray*}
S=\frac{|\nabla f|}{2} &=& \sqrt{x_n^2x_1^2+ \cdots+ x_n^2 x_{n-1}^2+ x_n^2 (x_n - k)^2 + ((x_n -k)k' + r r')^2}\\
&=& \sqrt{x_n^2 r^2 + A^2}.
\end{eqnarray*}
Thus we have
\begin{eqnarray*}
-nHS^3 = 2A^2x_n^2 +(n-1)kr^2x_n^3 +(n-2)kx_n A^2 + r^2x_n^2B - 2x_n^3k'A +2kk'x_n^2A .
\end{eqnarray*}
Squaring the above equation, we obtain
\begin{align} \label{eqn:mean curvature}
n^2 H^2 S^6 &= n^2 H^2 \{x_n^2 r^2 + (k' x_n + r r'- kk' )^2\}^3  \nonumber \\
&= \Big[\{2rr'k' + (n-2)kk'^2 +(n-1)kr^2 -r^2k''\}x_n^3  \\
& \quad +  \{(k'^2+r'^2-rr''+ kk'')r^2+2(n-3)kk'rr'-2(n-2)k'^2 k^2\}x_n^2 \nonumber \\
& \quad +  k(n-2)(rr'-kk')^2 x_n\Big]^2 \nonumber
\end{align}
Suppose that $H \neq 0$. Let us fix a section $t$. Since $x_n$ is varied,
we regard (\ref{eqn:mean curvature}) as an equation on $x_n$
where the coefficients are functions of the independent variable $t$.
Comparing the degree $0$-terms in both sides of (\ref{eqn:mean curvature}), we get
\begin{eqnarray*}
n^2H^2 (rr'-kk')^6 = 0.
\end{eqnarray*}
Therefore it follows that
\begin{eqnarray*}
\frac{dK(t)}{dt}=\frac{d}{dt} \sqrt{k(t)^2 - r(t)^2}=\frac{kk'-rr'}{\sqrt{k^2 - r^2}}=0.
\end{eqnarray*}
Now suppose that $H = 0$ and $n \geq 3$. Comparing the coefficients of the degree $2$-terms
in both sides of (\ref{eqn:mean curvature}), we have
\begin{eqnarray*}
k(n-2)(rr'-kk')^2=0.
\end{eqnarray*}
Therefore $rr'-kk'=0$, which also implies that
\begin{eqnarray*}
\frac{dK(t)}{dt}=0.
\end{eqnarray*}
Hence we can conclude that $M$ is a
rotationally symmetric hypersurface in both cases.

\qed \\

\begin{Rmk}
In $\mathbb{H}^2 \times \mathbb{R}$, Hauswirth \cite{Hauswirth06}
constructed several Riemann type minimal surfaces foliated by
circles. However, as mentioned in the introduction, it follows
from the above theorem that there is no Riemann type minimal
hypersurface which is not rotationally symmetric and foliated by
$(n-1)$-dimensional spheres lying in parallel horizontal
hyperplanes in $\mathbb{H}^n \times \mathbb{R}$ when $n\geq 3$.
\end{Rmk}
\section{Sphere-foliated hypersurfaces in the Lorentzian product $\mathbb{H}^n \times
\mathbb{R}$}

An immersed hypersurface $M$ in the Lorentz product space $\mathbb{H}^n \times
\mathbb{R}$ endowed with the Lorentzian metric
\begin{eqnarray*}
ds^2 = \frac{(dx_1)^2 + \cdots + (dx_{n})^2}{{x_n}^2} - (dt)^2
\end{eqnarray*}
is called $spacelike$ if the induced metric on $M$ is a Riemannian metric. If the
hypersurface is locally the level set of a smooth function $f$, the fact that $M$ is
spacelike means that $\nabla f$ is a $timelike$ vector:
\begin{eqnarray*}
\langle\nabla f, \nabla f\rangle < 0.
\end{eqnarray*}
If we orient $M$ by the unit normal vector field $\displaystyle{N=-\frac{\nabla f}{|\nabla f|}}$,
then the mean curvature $H$ is given by
\begin{eqnarray} \label{eqn:H in L}
nH = -\Div \frac{\nabla f}{|\nabla f|},
\end{eqnarray}
where $|\nabla f|=\sqrt{-\langle \nabla f, \nabla f \rangle}$ and $\Div$ denotes
the divergence with respect to the Lorentzian metric on the product space $\mathbb{H}^n \times
\mathbb{R}$.

As in the proof of Theorem \ref{thm:foliation}, consider two horizontal
hyperplanes of the foliation
$P_{t_1} = \mathbb{H}^n \times \{t_1\}$ and $P_{t_2} =
\mathbb{H}^n \times \{t_2\}$ for $t_1 < t_2$. Applying the Aleksandrov reflection
principle in ${\mathbb{H}}^n \times \mathbb{R}$,  we see that the
piece $M^*$ between $P_{t_1}$ and $P_{t_2}$ has the symmetries of its boundary
 $\partial M^* = (M^*\cap P_{t_1}) \cup (M^*\cap P_{t_2})$. Therefore, for each $t_1 \leq t \leq t_2$, the hyperbolic center
 of each level $M\cap \{x_{n+1}=t\}$ lies in the same $2$-plane. After a translation
in ${\mathbb{H}}^n \times \mathbb{R}$, we may assume that this $2$-plane is defined by $x_1= \cdots=x_{n-1}=0$.

Using Lemma \ref{lem:center} again, we can parametrize the hyperbolic centers of the
$(n-1)$-dimensional spheres by $t \longmapsto (0, \ldots, 0, K(t),
t) \in \mathbb{H}^n \times
\mathbb{R}$ for $t\in [t_1, t_2]$, and the Euclidean centers of the
spheres by $t \longmapsto (0, \ldots, 0, k(t), t)\in \mathbb{H}^n \times
\mathbb{R}$, where $K(t) = \sqrt{k(t)^2 - r(t)^2}$ and $r(t)$ is the
Euclidean radius of $M \cap \mathbb{H}^n \times
\{t\}$. Then $M^*$ is the level set $\{f = 0\}$ of a function $f$ defined as in (\ref{eqn:function}).
Note that the metric on the Lorentzian product $\mathbb{H}^n  \times \mathbb{R}$ is given by
\begin{eqnarray*}
\sum_{i,j} g_{ij} dx_i \otimes dx_j = \frac{1}{x_n^2}dx_1^2 + \cdots + \frac{1}{x_n^2}dx_n^2 - dt^2 .
\end{eqnarray*}
Since
\begin{eqnarray*}
\nabla f &=& \sum_{i,j} g^{ij}\frac{\partial f}{\partial x_i}\frac{\partial}{\partial x_j} \\
 &=& \Big(2x_n^2x_1, \cdots, 2x_n^2 x_{n-1}, 2x_n^2 (x_n - k), 2((x_n -k)k' + r r')\Big) ,
\end{eqnarray*}
we get
\begin{eqnarray*}
-\langle \nabla f, \nabla f \rangle =  4\Big(-x_n^2x_1^2- \cdots - x_n^2 x_{n-1}^2 - x_n^2 (x_n - k)^2 + ((x_n -k)k' + r r')^2\Big) .
\end{eqnarray*}
Using the mean curvature equation (\ref{eqn:H in L}), we have
\begin{eqnarray*}
-nH = \Div \frac{\nabla f}{|\nabla f|} &=& \sum_{i,j} \frac{1}{\sqrt{|g|}} \frac{\partial}{\partial x_j}\Big(\sqrt{|g|}\frac{g^{ij}}{|\nabla f|}\frac{\partial f}{\partial x_i}\Big) ,
\end{eqnarray*}
where $|\nabla f| = \sqrt{-\langle \nabla f, \nabla f \rangle}$ and $|g|=|\det(g_{ij})|=x_n^{-2n}$.
\vspace{3mm}

A similar computation as in the proof of Theorem \ref{thm:foliation} shows that
\begin{align} \label{eqn:mean curvature in L}
n^2 H^2 &\{-x_n^2 r^2 + (k' x_n + r r'- kk' )^2\}^3  \nonumber \\
&= \Big[\{2rr'k' - (n-2)kk'^2 -(n-1)kr^2 + r^2k''\}x_n^3  \\
& \quad +  \{(k'^2+r'^2-rr''+ kk'')r^2-2(n-1)kk'rr'+2(n-2)k'^2 k^2\}x_n^2 \nonumber \\
& \quad -  k(n-2)(rr'-kk')^2 x_n\Big]^2 \nonumber
\end{align}
Suppose that $H \neq 0$. Comparing the degree $0$-terms in both sides
of the equation (\ref{eqn:mean curvature in L}) of variable $x_n$, we obtain
\begin{eqnarray*}
n^2H^2 (rr'-kk')^6 = 0.
\end{eqnarray*}
Therefore it follows that
\begin{eqnarray} \label{eqn:center01}
\frac{dK(t)}{dt}=\frac{d}{dt} \sqrt{k(t)^2 - r(t)^2}=\frac{kk'-rr'}{\sqrt{k^2 - r^2}}=0.
\end{eqnarray}
Now suppose that $H = 0$ and $n \geq 3$. Comparing the coefficients of the degree $2$-terms
in both sides of (\ref{eqn:mean curvature in L}), we have
\begin{eqnarray*}
k(n-2)(rr'-kk')^2=0.
\end{eqnarray*}
So we have $rr'-kk'=0$, which also implies that
\begin{eqnarray} \label{eqn:center02}
\frac{dK(t)}{dt}=0.
\end{eqnarray}
From (\ref{eqn:center01}) and (\ref{eqn:center02}), it follows that the hyperbolic
center of each hypersphere in parallel horizontal hyperplane lies in a vertical
geodesic line of the Lorentzian product $\mathbb{H}^n \times \mathbb{R}$.
Therefore we obtain the following.
\begin{thm}
Let $M$ be an $n$-dimensional spacelike hypersurface with constant
mean curvature $H$ in the Lorentzian product $\mathbb{H}^n \times
\mathbb{R}$ and foliated by spheres in parallel horizontal
hyperplanes. If $H \neq 0$ or $H=0$ and $n \geq 3$, then $M$ is
rotationally symmetric.
\end{thm}


\vspace{1cm}
\ \  \\
Department of Mathematics\\
Sookmyung Women's University\\
Hyochangwongil 52, Yongsan-ku \\
Seoul, 140-742, Korea\\
e-mail : kseo@sookmyung.ac.kr
\end{document}